\documentclass{svmult}

\usepackage{mathptmx}       
\usepackage{helvet}         
\usepackage{courier}        
\usepackage{type1cm}        

\usepackage{graphicx}        
\usepackage{multicol}        
\usepackage[bottom]{footmisc}

\usepackage{amsfonts,amsmath,latexsym,amssymb,makeidx,multicol,supertabular}
\usepackage[colorlinks=true,linkcolor=black,citecolor=black,urlcolor=black]{hyperref}

\usepackage[top=1.4cm, bottom=2.4cm, left=1.9cm, right=1.5cm]{geometry}
\def\:{:\,}
\def\definedby{\stackrel{\Delta}{=}}

\def\nomgroup#1{}
\def\indic{{{1}}}

\def\EC{Euler characteristic\ }

\def\bb{\mathbb }

\def\cal{\mathcal}

\def\det{{\mbox{\rm det}}}

\def\LKC{Lipschitz-Killing curvature}

\def\P{{\bb P}}
\def\E{{\bb E}}

\def\I{{\rm Ind}}

\def\p{\varphi} 
\def\real{{\bb{R}}}

\def\varPsi{\Psi}
\def\RN{\mathbb{R}^N}

\def\Rk{\mathbb{R}^k}

\def\definedas{\stackrel{\Delta}{=}}

\def\Min{{\cal M}}

\newcommand{\sqbinom}[2]{\begin{bmatrix}#1 \\ #2 \end{bmatrix}}

\def\tt{\widetilde t}

\def\cO{{\cal O}}

\newcommand{\bc}{\begin{center}}
\newcommand{\ec}{\end{center}}

\newcommand{\eex}{\end{exercise}}
\newcommand{\bex}{\begin{exercise}}
\newcommand{\beq}{\begin{eqnarray}}
\newcommand{\eeq}{\end{eqnarray}}
\newcommand{\beqq}{\begin{eqnarray*}}
\newcommand{\eeqq}{\end{eqnarray*}}

\newcommand{\mink}{{\cal M}}

\newcommand{\lips}{{\cal L}}

\newcommand{\Tube}{{\rm Tube}}

\def\text{\mbox}

\newcommand\Var{\mathop{\rm Var}}

\begin{document}

\title*{Some New Random Field Tools for Spatial Analysis}
\titlerunning{Random Fields} 
\author{Robert J.\ Adler}
\institute{Robert J.\ Adler \at Faculties of Industrial Engineering \& Management and of  Electrical Engineering, Technion, Haifa, Israel, 32000\\ \email{robert@ieadler.technion.ac.il}}

\maketitle
\begin{multicols}{2}
\abstract{This is a brief review, in relatively non-technical terms, of recent 
rather technical advances in
the theory of random field geometry. These advances have provided a collection of explicit
new formulae describing mean values of a variety of geometric characteristics of 
excursion sets of random fields,  
such as their volume, surface area and Euler  characteristics.
What is particularly important in these formulae is that whereas the previous theory
covered only stationary, Gaussian random fields, the new theory 
requires neither stationarity nor, a fortiori, isotropy. Furthermore, it
covers a wide class of non-Gaussian random fields.
The  formulae provided by these advances have a wide range of potential applications, 
including new techniques of 
parameter estimation, model testing, and thresholding for spatial and space-time 
functions.
The paper reviews this theory, and provides brief descriptions of some of the
 applications.}


\vskip0.6cm \hrule height  0.5pt
\section{Introduction}
\label{sec:intro}

Random fields, a generic term that I shall use 
to describe random processes defined over parameter
spaces with dimension great than one, appear throughout the modelling of spatial and 
space-time phenomenoma, whether space is two or three dimensional. To note just
a handful of examples, we have
\begin{itemize}
\item Physical oceanography and hydrology, where,  as in \cite{AMR,Rubin},
 the spatial parameter might be the 
two dimensional water surface or a three dimensional body of water
 and the variable being measured might be water temperature and/or
pressure. If both are being measured, then we have an example of a vector-valued random
field.
\item Atmospheric studies, where the random field might provide a model for 
wind speads or airborne contaminants, with a theory as described in 
\cite{Daley}.
\item Geostatistics and other earth sciences, where, as for oceanography, the 
spatial aspects of parameter spaces may be two or three dimensional. cf.\
\cite{Christakos1,Christakos2,Olea}.
\item Astrophysics, where random field techniques have been heavily used in analysing the 
COBE (\underline{Co}smic \underline{B}ackground \underline{E}xplorer) data, which 
measures the `signature radiation' from the universe of 15 billion years ago
\cite{Smootbook,TORRES,VPGHG}. This is directional data, and so is realised as a
random field on the two-dimensional sphere. Three dimensional astrophysical
data has been generated by the Sloan digital sky survey   \cite{Sloan} where, as
we shall describe later, heavy use has been made of the geometrical theory that is
at the centre of this review.
\item Analysis of functional magnetic resonance imaging (fMRI)  data. This is the archtypical example that will be used 
throughout this paper, and so I shall leave the details until later.
\end{itemize}

The theory of random fields has undergone a minor revolution over the past few years in terms
of its ability to provide precise geometric information about the global behaviour
of sample paths. While
this work was motivated by applications, it has appeared primarily in the pure
mathematical literature and, even by its esoteric standards, is not always very easy to read.

In this brief review we want to describe some of these results in simpler language,
 as well as discussing some of their applications.

The new results centre on the geometry of the so-called 
{\it excursion sets} or  {\it nodal domains} of random fields. To describe these, we first
define a {\it random field} to be a random function $f$ defined over a parameter space $M$, which
we shall always take to be a bounded region in $N$-dimensional Euclidean space, $\RN$.
The random field itself might be real or vector valued. In general, we let it take values in
$\Rk$, $k\geq 1$. If $D$ is a subset of $\Rk$, then the excursion set of $f$ in $D$ and over
$M$ is defined by
\beq
\label{intro:AD:equn}
A_D \equiv A_D(f,M) \definedas \{x\in M\: f(x)\in D\}.
\eeq
In the particular case that $f$ is real-valued and $D$ is the set $[u,\infty)$, we are looking at
those points in the parameter space $M$ at which $f$ takes values larger than $u$, and write
\beq
\label{intro:Au:equn}
A_u \equiv A_u(f,M) \definedas \{x\in M\: f(x)\geq u\}.
\eeq
In general, we shall call $D$ a {\it hitting set}.

Perhaps the first studied example of a real valued random field in two dimensions 
is due to Mat\'ern \cite{MATERN}, who modelled forestry yields $f(x)$ 
as a function of a two-dimensional positional variable $x$. The excursion sets were then
high yield geographical regions. Real valued 
random fields defined over $\real^3$ abound in environmental
data, given, for example, by pollutant levels in space. If more than one pollutant is 
measured, then the random field becomes vector valued and the sets $D$ defining
excursion sets can become quite involved. For example, if interest is in regions where at 
least one of the pollutants takes some minimal value -- say $u_i$ for the $i$-th pollutant,
$i=1,\dots,k$ -- then $D$ will be of the rather simple, rectangular 
form $\prod_{i=1}^k [u_i,\infty)$. However, if interest lies in some combination of pollutants,
then $D$ will reflect this and will be, geometrically, as complicated as is the formula for
the combination.

Examples with four dimensional parameter sets are also common: One needs only add time to 
examples that otherwise deal with three dimensional space. Space-time examples are important, and
an example where the new theory has something to offer that the old theory did not. It is
obvious that, in essentially all real scenarios, the structure of temporal dependence is quite
different to that of spatial dependence. Therefore, models restricted to assumptions such as
isotropy become problematic in space-time. Indeed, even stationarity can be a problem in
dealing with processes where, perhaps because of external factors,  there is no intrinsic 
temporal stationarity, despite there being spatial stationarity, or vice versa. 
We shall close the introduction with a (perhaps surprisingly) simple
example of an eight-dimensional random field, which highlights the need for a theory that
does away with classical assumptions on the simplicity of the parameter space and 
directionally homogeneous stochastic behaviour of the random field.

However, since we have not described what the new theory does, we shall do this first.
Given that excursion sets are of intrinsic interest, it would be nice to understand a little of
their properties, such as how many components they have, how large they are, what is the size
of their boundaries, etc. There are a number of ways to quantify these properties, and for $N$
dimensional sets geometers have developed 
$N+1$ numerical quantifies, known under a variety of names, 
including  Quermassintegrales, Minkowski or Steiner functionals, integral curvatures,
 intrinsic volumes, and {\em \LKC s}, hereafter LKCs. The differences between them are only of 
ordering and scaling, and
although it is the longest of all names, we shall choose LKCs, for consistency with our basic reference \cite{RFG}.
 We shall denote the LKCs of a $N$-dimensional set $A$ by $\lips_0(A),\dots,\lips_N(A)$ and, in many ways, they can  be thought of as measures of 
the  `$j$-dimensional sizes'
of $A$.  For example, when
$N=2$ they have the following meanings:
\begin{itemize}
\item $\lips_2(A)$  is the two dimensional area of $A$.
\item $\lips_1(A)$  is half the boundary length of $A$.
\item $\lips_0(A)$ is the \EC\ of $A$, which in two dimensions is given by
\beq
\label{EC2}
\lips_0(A) &=& \#\{\text{connected components in $A$}\} \qquad \\ \notag
&&\qquad\qquad - \ \#\{\text{`holes' in $A$}\}.
\eeq 
\end{itemize}
When $N=3$,
\begin{itemize}
\item $\lips_3(A)$  is the three dimensional volume of $A$.
\item $\lips_2(A)$  is half the surface area of $A$.
\item $\lips_1(A)$  is twice the {\it caliper diameter} of $A$,
where the caliper diameter of a convex $A$ is defined by placing the solid
between two parallel planes (or calipers), measuring the distance between the
planes, and  averaging over all rotations of $A$. 
Thus the caliper diameter of a sphere is just the
usual diameter, while for a rectangular box of size $a\times b\times c$ 
it is $(a+b+c)/2$, half
the `volume' used by airlines to measure the size of luggage.
\item $\lips_0(A)$ is again the \EC\ of $A$, which in three dimensions is given by
\beq
\label{EC3}
\lips_0(A) &=& \#\{\text{connected components in $A$}\}\qquad \quad \\ \notag
&& - \#\{\text{`handles' in $A$}\} + \#\{\text{`holes' in $A$}\}.
\eeq 
\end{itemize}
We shall give further definitions for general dimensions in Section \ref{sec:geometry}.

What the new theory provides is explicit formulae for the expectations of the LKC's of
excursion sets -- that is for
\beq
\label{ELKC.1}
\E\left\{\lips_j\left(A_D(f,M)\right)\right\},
\eeq
--  for a wide class of random fields $f$, over a large class of parameter sets $M$, and
for a large class of sets $D$.

Without yet being explicit about the precise structure of these formulae, we discuss some of
their uses.

\subparagraph{Signal detection and thresholding}
Consider the classical  signal+noise paradigm, which  considers the model
\beq
\label{signal-noise}
f(x) = s(x) + \eta (x).
\eeq
The signal $s$ is  deterministic, and, purely for reasons of exposition,
is assumed to  take primarily 
positive values if, in fact, it is present. The noise $\eta$ is a mean zero 
random field of background noise, the statistical properties of which are assumed to be known. 
The first problem here is determining whether or not the signal is at all present, and,
in view of the assumed positivity of $s$, a
useful statistic for testing this is the threshold statistic
\beqq
\sup_{x\in M} f(x).
\eeqq
 (In fact, the supremum is the maximum likelihood statistics for this
test, if  $\eta$ is Gaussian white noise smoothed with a filter with shape matching that
of the signal. This is the `matched filter theorem' of signal processing.)

In order to use this statistic, or, indeed, to perform any statistical thresholding technique, one needs to know how to compute, at least for large values of $u$, the {\em excursion probabilities}
\beq
\label{excursionprob}
\P\Big\{\sup_{x\in M} f(x)\, \geq\, u\Big\},
\eeq
when no signal is present.
(By `large' we mean large enough for this probability to be less that 
about 0.10, the level at which one usually begins to talk about statistical significance.)

It has long been acknowledged that, for many random fields,  
\beq
\label{ECequivalence}
&& \Big|\P\Big\{\sup_{x\in M} f(x) \geq u\Big\} - 
\E\left\{\lips_0\left(A_u(f,M)\right)\right\}\Big|
\leq error(u),
\notag \\
\eeq
where the $error$ function is of a smaller order than both of the other terms. This means that
the excursion probabilities \eqref{excursionprob}, which are intrisically uncomputable, can
effectively be replaced by an expectation which might be more accessible.

Two  of the contributions of the recent theory are the computation of this expectation
 in wide generality, and a
full proof that, at least for Gaussian random fields, not only finally rigorously established
 \eqref{ECequivalence} but even identified
 the structure of the function $error$.
(cf.\ \eqref{errorf} below.)

\subparagraph{Parameter estimation}

Assume that we are dealing with a random field, the general  distributional structure of which is 
known, but for which values of  parameters are unknown. There are standard statistical techniques for
estimating the means and variances of  random fields, as well as for  estimating covariance
functions, which we shall denote as
\beqq
C(x,y) = \E\{f(x)f(y)\},
\eeqq
assuming, from now on, and for notational convenience, that any  non-zero mean has already
by removed from $f$. (When $f$ is stationary, we shall adopt the usual notational
inconsistency   of writing $C(x-y)\equiv C(x,y)$.)
However, estimates of covariance functions are often not well behaved,
particularly in the non-stationary case, when they are  notoriously unreliable. Furthermore,
in many situations, one needs to know only certain aspects of the covariance function, and
not the entire function itself. For example, if one is concerned with excursion probabilities
\eqref{excursionprob} for large $u$, and if $f$ is stationary and Gaussian, then it turns
out that one needs only to know certain derivatives of $C$ in the neighbourhood of the 
origin, or, equivalently, the second spectral moments of $f$. 
(cf.\  \eqref{lambda2}, \eqref{lambdaij} below.)

As opposed to simple means and variances, which describe static properties of $f$ (i.e.\
properties of $f$ at specific points) covariance and 
 spectral parameters describe dynamic behaviour over the entire parameter space, and much can be 
done with them, without knowing the full covariance function. From a statistical point of view,
it is, of course, more efficient to estimate a handful of parameters than an entire function.

However, it is precisely these parameters that appear in the formulae for the expected LKCs of 
excursion sets, and so those formulae provide a tool for  estimating these parameters, by
comparing expectations to empirically observed properties of excursion sets.

\subparagraph{Model testing}

In this problem, which is in some sense an extension of the parameter estimation problem,
we are in the scenario of dealing with a random field whose precise structure is unknown, but 
it assumed to belong to some small class of models. Since the forms of the expected LKCs of
excursion sets are often quite different for different models, comparing empirical data
to various theoretical models can be used to choose among models. 
This technique has,
for example, been used quite heavily in the astrophysical literature (e.g.\
\cite{Smootbook,TORRES,VPGHG}) and we shall give an example of how to use it in Section 
\ref{sec:examples}  below.

Finally, we complete this introduction of a quite concrete example of a random field over
a rather complicated subset of $\real^8$, in an attempt to convince the possibly sceptical, 
applied reader that  there is actually is a need for the practioners' tool kit
to include random fields over high dimensional parameter spaces.

\subparagraph{An eight dimensional example in scale-space} We shall take an example from
fMRI brain imaging, since we can then refer you to \cite{SSSW} for further details, where
the approach is notationally quite consistent with this review. The problem, however, is 
quite generic.

In fMRI procedures, measurements are taken at a 
large number of voxels in the (three-dimensional) 
brain. For simplicity, we
shall assume one such measurement per voxel. However, these are not actually point
 measurements, but each measurement is rather a local average, taken over a small 
ellipsoidal neighbourhood
 of the voxel. Each such ellipsoid has three principle axes, the lengths of each of 
which can be chosen by the technician operating the scanner, giving three more dimensions to
the problem. The final two dimensions come from possible rotations of the ellipsoid.
Thus, what seemed originally to be an observation taken at a point in the three dimensional
brain, is really an observation taken in the eight dimensional space 
\beq
\label{brain}
\text{brain\,$\times$\,axes\,$\times$\,rotation.}
\eeq

Why is this important? Suppose we are interested in either  a signal detection or 
thresholding problem as described
above, looking for signs of unusual behaviour somewhere in the brain. The technician, as he
plays with the sizes of the axes and their
 rotations, is clearly playing with the resolution of the image, and
the classic `uncertainty principle' comes into play. The higher the resolution (i.e.\ the
smaller the lengths of the principle axes of the ellipsoids) the lower is the statistical
reliability of  the measurements at individual voxels. Thus, in changing these parameters, the technician is 
implementing a control procedure which can impact quite strongly on the excursion 
probability \eqref{excursionprob} and so on any thresholding and risk assessment
 that comes from it, unless this probability
 is calculated in the full, eight dimensional, scenario.

Note that the covariance function of this random field will be neither stationary nor
isotropic. Even if one treats the basic, three-dimensional 
 structure of the brain as isotropic (and even this is suspect) the additional 
dimensions added here behave very differently. Consequently, handling this example calls
for the development of a theory of non-stationary random fields.

Note also that while we chose the fMRI example to make our point, the problem itself is
generic to any situation in which data collection involves a compromise between 
resolution and smoothness. To put things more bluntly, the problem arises whenever raw
spatial data is made up of some sort of local average, which is an almost ubiquitous
situation.

\vskip0.6cm \hrule height  0.5pt
\section{Gaussian and Related Random Fields}
\label{sec:rfields}
The basic random fields with which we shall work are Gaussian ones, so that both their univariate
{\em and multivariate} distributions are all normal. In the real valued case, this means that 
the  random function $f\:M\to\real$  has the property that all collections
\beqq
          f(x_1),\dots, f(x_k),
\eeqq
of random variables have multivariate Gaussian distributions, for all $k\geq 1$ and all
 $x_1,\dots,x_k\in M$. A similar definition holds for vector valued random fields.
 We have emphasised the words `and multivariate' above
on purpose. Throughout the literature one finds examples in which authors make a pointwise
transformation to normality, and then claim normality for their process. The argument, for
real valued random fields, usually goes as follows:

Collect data points of the random field $f$ at points $x_i\in M$, $i=1,\dots,n$, 
and form the   empirical distribution function
\beqq
\widehat F_n(u) = \frac{\#\{i\: f(x_i)\leq u\}}{n}.
\eeqq
Letting  $\Phi$ denote the distribution function of a standard normal 
random variable, define the  `Gaussianised', or sometimes just 
`standardised', data to be
\beq
\label{intro:standard:equn}\qquad
\widetilde f(x_j) \definedby \Phi^{-1}\big(\widehat F_n(f(x_j))\big),
\quad j=1,\dots,n.\qquad
\eeq
If the random field is stationary, and the sample size is large enough, then it certainly 
is true that the random field $\widetilde f$   has univariate distributions that are close to
standard normal. However, no claim at all can be made about the {\it multivariate}
 distributions, 
and so, a fortiori,  about $\widetilde f$ as a {\em process}, and the application of Gaussian 
process theory can lead to quite significant errors. (The same is true even if the empirical
distribution function in \eqref{intro:standard:equn} is replaced with the true, theoretical
one, as we shall see in Section  \ref{sec:examples}.)

\subsection{Gaussian related fields}
\label{sec:grelated}
Leaving the Gaussian scenario is, however, not all that easy to do, and the new theory
leaves it in a fashion that, while somewhat limited, turns out to be  broad enough
to cover many, if not most, statistical applications of random fields. To be more precise,
we shall call a random field $f:M\to \real^d$ a {\it Gaussian related}
 field if we can find 
a vector  valued Gaussian random field,
\beqq
g(x)=\left(g_1(x),\dots,g_k(x)\right)\: M\to\real^k,
\eeqq
 with independent, identically distributed components having zero means and unit variances, and a function
\beq
\label{F}
F\:\real^k\to\real^d,
\eeq
such that $f$ has the same multivariate distributions as $F(g)$.

When $k=1$, or, in general $k=d$ and $F$ is invertible, then the corresponding
Gaussian related process is not much harder to study than the original Gaussian
one, since  what happens at the level $u$ for $f$ is precisely what happens at
the uniquely defined level $F^{-1}(u)$ for $g$. In fact, this 
actually falls into the class of transformations that we warned against above,
but in this case the transformation argument ifs completely valid.

In the more interesting cases in which  $F$ is not invertible, $f=F(g)$ can
provide a process that is qualitatively different to $g$.
 For example,  consider 
the following three choices for $F$, where in the third we set $k=n+m$.
\beq
\label{nongauss:threeexamples:equation}
\sum_1^k x_i^2,\qquad
\frac{x_1\sqrt{k-1}}{(\sum_2^k  x_i^2)^{1/2}},\qquad
 \frac{m\sum_1^n  x_i^2}{n\sum_{n+1}^{n+m}  x_i^2}.
\eeq
The corresponding random fields are  known as  $\chi^2$ fields with 
$k$ degrees of freedom, the $T$ field with 
$k-1$ degrees of freedom, and the  $F$
  field with $n$ and $m$ degrees of freedom. These three random fields all have very different
spatial behaviour, and each is 
as fundamental to the statistical applications of random field theory as is its corresponding
univariate distribution to standard statistical theory. In each of these three cases, as
in general for a Gaussian related random field, there is no simple pointwise transformation
which will transform   it to a real valued Gaussian field.

Note that for a Gaussian related field $f$ the excursion sets $A_D$ 
can be rewritten as
\beq
\label{intro:nongaussexc:equation}
 A_D(f,M) &=& A_D(F(g),M)\\ &=& \{x\in M\: (F\circ g)(t) \in D\}   \notag  \\
&=& \{x\in M\: g(x) \in F^{-1}(D)\} \nonumber  \\
&=&  A_{  F^{-1}(D)}(g,M).        \nonumber 
\eeq
Thus, for example, the excursion set of a {\it real valued}
  non-Gaussian $f=F\circ g$ above a level
$u$ is equivalent to the excursion set for a {\it vector valued} 
 {Gaussian} $g$ in $F^{-1}([u,\infty))\in \Rk$. 

We shall return to this important point later.

\subsection{Regularity Assumptions}
\label{regularity:sec}
Although we do not want to get into detailed technicalities here, we do need
to put some restrictions on the central objects of this review. We start
with the parameter spaces.

\subparagraph{Parameter spaces}
The general theory allows $M$ to be what is known as a 
 Whitney stratified manifold, satisfying some mild side conditions. Rather
than define what these are, we shall give a list of examples which should
 cover most of the situations found in environmental engineering.
\begin{itemize}
\item[(i)] Finite simplicial complexes. This includes sets such as 
$N$-dimensional rectangles of the form 
\beq
\label{rectangles}
T=\prod_{i=1}^M [0,T_i],
\eeq
or unions of a finite number of such rectangles,
as well as most sets with flat boundaries.
\item[(ii)]\ Reimannian manifolds, with or without boundary. This includes
sets such as $N$-dimensional balls and spheres, and smooth deformations of them.
\item[(iii)]\ \ \ Sets that can be written in the form
\beq
\label{stratification}
M=\bigcup_{i=1}^N \partial_iM,
\eeq
where $\partial_i$ is an $i$-dimensional set (the $i$-dimensional `boundary'
of $M$) that fits into the one of 
the first two categories. For example, the eight dimensional example of
\eqref{brain}, which fits into neither of the first two categories, 
 requires this level of generality.
\end{itemize}

We place on the hitting sets $D$ appearing in the definition
\eqref{intro:AD:equn} of excursion sets one out of the above possible 
 assumptions on $M$.

\subparagraph{Random fields}
The basic building blocks of all the random fields we consider are smooth
Gaussian random fields, $f\:M\in\RN\to\real^k$. The first assumption is a
minor one, and only for notational convenience, 
that all means be fixed at zero. The second is more serious, and it is
that $f$ have constant variance throughout $M$. Note that this is a much
weaker assumption than stationarity, and is achievable in general by 
replacing a random field that does not have constant variance by
\beqq
f(x)\to f_{standardised}(x) = \frac{ f(x)}{(\E\{f^2(x)\})^{1/2}}.
\eeqq
This transformation is often useful and involves no loss of infomation. For example,
the excursion sets above zero are the same for $f$ and $f_{standardised}$.

However, the fact that the new theory requires only constant variance and not
full stationarity is extremely useful, since in general there is 
no simple transformation which will take a non-stationary field to a stationarity
one. (See, however, \cite{SampsonGuttorp}, which
shows that if one is prepared to pay the price of moving to a higher
dimensional scenario than the one inherent to the problem, there is a 
way around this.) 

The final assumption of consequence on the Gaussian processes, and the most 
important one,
 is that each of the $k$ components of
$f$ is twice differentiable, and that these derivatives are themselves
continuous. Some additional minor assumptions of non-degeneracy also
need to be made, but since these are of little practical importance we
direct the interested reader to Chapter 11 of \cite{RFG} or Chapter 4 of \cite{ARF} for details.

\subparagraph{Gaussian related fields}
Recall that Gaussian related fields were defined as pointwise transformations
of vector valued Gaussian fields. We have already assumed that the components
 of these fields are independent and identically distributed, and now add the assumption
that  the  
function $F$ of $\eqref{F}$ be twice continuously differentiable. As 
a consequence, all our Gaussian related fields are also 
 twice continuously differentiable.

\vskip0.6cm \hrule height  0.5pt
\section{Tube formulae}
\label{sec:geometry}
In the Introduction we described the meanings of the \LKC s $\lips_j$ for two and 
three dimensional sets, but we also need a definition in higher dimensions as
well. We shall also soon need another set of related geometric quantifiers, 
which we call Gaussian Minkowski functionals, and all of these will be defined
in this section.
 
\subsection{\LKC s}
\begin{figure}[!ht]
\rotatebox{-90}{\includegraphics[scale=.25]{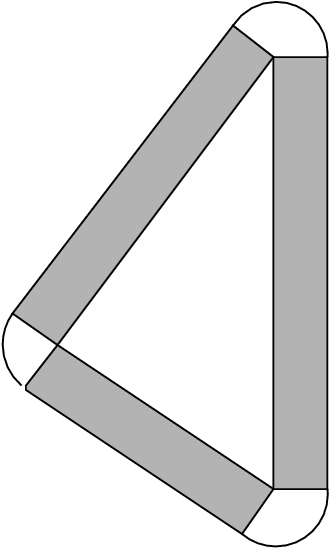}}
\caption{The tube around a triangle.}
\label{figtriangle}
\end{figure}

The quickest way to define \LKC s is via a basic result of integral geometry 
known as {\em Steiner's formula}. Steiner's formula deals with the $N$-dimensional
volume of enlargements of sets, where the enlargment, or {\it tube}, of diameter
$\rho > 0$  built around a set $A\in\RN$ is the set
\beq
\label{steiner}
\Tube(A,\rho) = \Big\{x\in \RN\: \min_{y\in A} \|x-y\|\leq \rho\Big\}. \quad 
\eeq
An example is given in Figure \ref{figtriangle}, in which $A$ is a
solid  triangle and $\Tube(A,\rho)$ is the larger triangular region with
rounded corners. Note that $\Tube(A,\rho)$ always includes $A$ itself. To see
from where the terminology comes, think of $A$ as being a curve in $\real^3$,
in which case $\Tube(A,\rho)$ really does look like a (solid) physical tube.
\begin{center}
\begin{tabular}{|c|}   \hline   \\
{\ \ Place Figure \ref{figtriangle} near here.\ \  }
\\ \\   \hline \end{tabular}
\end{center}
Steiner's theorem says that if $\lambda_N$ denotes volume in $\RN$, then, for
convex $A$ of dimension $\dim (A)$,
\beq
\label{steiner}
\lambda_N
\left(\Tube (A,\rho)\right)  = \sum_{j=0}^{\dim (A)} \omega_{N-j} \rho^{N-j} \lips_j(A), \quad
\eeq
where $\omega_j ={\pi^{j/2}}/{\Gamma\bigl(1+ {j}/{2} )}$
is the volume of the unit ball in $\real^j$.
This (finite) power series in $\rho$ defines the LKCs, and by looking at 
examples, such as squares in the plane and cubes in $\real^3$, you should find
is easy to convince yourself that the descriptions given earlier of the LKCs in
two and three dimensions are reasonable.

There are versions of Steiner's formula for far more general sets than convex
ones, in which case they are usually called `tube formulae' and generally associated
with the name of Hermann Weyl. For us it will suffice to know that, for small
enough $\rho$, \eqref{steiner} holds for all the sets that we shall consider in
this review.
 
It is quite easy to see from \eqref{steiner} that the LKCs satisfy a basic
scaling relationship, in that 
\beqq
\lips_j(\lambda A) = \lambda^j\lips_j (A),
\eeqq
for $0\leq j\leq \dim (A)$ and $\lambda>0$, where 
$\lambda A= \{x\:x=\lambda y,\ \text{for some}\ y\in A\}$. Thus one can think
of $\lips_j(A)$ as some sort of measure of the `$j$-dimensional volume' of $A$.

\subsection{Gaussian Minkowski functionals}
\label{gmf:sec}
We need one more set of geometrical tools, closely related to LKCs, before we can describe
the new results that are at the core of this review, and these are related to the hitting
sets $D$ of \eqref{intro:AD:equn}.

Rather than measuring the basic Euclidean sizes of the hitting sets, we want to measure
aspects of their (Gaussian) probability content. In particular, let $X$ be a vector of $k$ 
independent, identically distributed, standard Gaussian random variables, and,
for a nice set $A\in\Rk$, write
\beqq
\gamma_k(A)=\P\{X\in A\}=\int_A \ \frac {e^{-\|x\|^2/2}}{(2\pi)^{k/2}}\, dx.
\eeqq
Defining tubes exactly as before, it turns out that there is a Taylor expansion for the 
probability contact of tubes (due to Jonathan Taylor \cite{Taylor1}), which we shall write as
\beq
\label{gaussexp1}
\gamma_{k}(\Tube(A, \rho))  =   \sum_{j=0}^{\infty} 
\frac{\rho^j}{j!} \Min^{k}_j(A),
\eeq
for small enough $\rho$.
We call the coefficients $ \Min^{k}_j(A)$ the {\em Gaussian Minkowski functionals}  (GMFs) of $A$.
Clearly $\Min_0^k(A)=\gamma_{k}(A)$. Note that, unlike the expansion \eqref{steiner}
in Steiner's theorem, this expansion is an infinite one.

To see how this expansion works in the simplest of cases, take $k=1$ and 
$A=[u,\infty)$.  Then, since 
\beqq
\gamma_1\big(\Tube([u,\infty),\rho)\big)=\gamma_1\big([u-\rho,\infty))\big)=\Psi(u-\rho),
\eeqq
where 
\beqq
\Psi(u)=(2\pi)^{-1/2}\int_u^\infty e^{-x^2/2}\,dx
\eeqq
 is  the tail probability function  for  a standard Gaussian  variable, a standard Taylor
expansion along with properties of Hermite polynomials shows that 
\beq
\label{tubes:minkequalshermite:equation}
\Min^{k}_j([u,\infty))
=  H_{j-1}(u) \frac{e^{-u^2/2}}{\sqrt{2\pi}}.
\eeq
Here $H_n$, $n\geq 0$, is the $n$-th Hermite polynomial
\beq
\label{rgeometry:Hermite:equation} \qquad
H_n(x) =  n!\, \sum_{j=0}^{\lfloor n/2\rfloor} 
\frac{(-1)^j x^{n-2j}}{j!\, (n-2j)!\, 2^j},
\eeq
where $\lfloor a\rfloor$ is the largest integer  less than or equal to
$a$. When $n=-1$ we set
\beq
\label{rgeometry:Hermite-zero:equation}
H_{-1}(x) =  \sqrt{2\pi}\, \Psi(x)\,e^{x^2/2}.
\eeq

In fact, it is not much harder to define the GMFs of more general sets that arise in the study
of Gaussian related random fields. With the notation of Section  \ref{sec:grelated}
let $f=F(g)$ be such a field. Then, as we already noted above, it follows from
 \eqref{intro:nongaussexc:equation},
the excursion sets $A_u(f,M)$ are those of $g$ in the hitting sets $F^{-1}([u,\infty))$, and
so these are going to be of interest to us. If we now write $p_F(u)$ for the probability 
density of the random variable $F(g(x))$ (which for simplicity we shall assume 
does not depend on $x$) then 
\beqq
\mink^k_j\big(F^{-1}([u,\infty))\big) = (-1)^{j-1}\frac{d^{j-1}p_F(y)}{dy^{j-1}}\Big|_{y=u},
\eeqq
for $j\geq 1$, while
\beqq
\mink^k_0\big(F^{-1}([u,\infty))\big) = \gamma_k\big(F^{-1}([u,\infty))\big).
\eeqq
Although these formulae seem simple, they often involve considerable delicate albegra to
compute. For example, if $F(x)=\|x\|^2$, so that the Gaussian related field is $\chi^2$ with
$k$ degrees of freedom, then, for $j\geq 1$, 
\beq
\label{chisqared}
&&\Min_j^k \big(F^{-1}([u,\infty))\big)\\&& =   \frac{u^{k-j}e^{-u^2/2}}{\Gamma(k/2) 2^{(k-2)/2}} \sum_{l=0}^{\lfloor \frac{j-1}{2} \rfloor} \sum_{m=0}^{j-1-2l}  
 \indic_{\{k \geq j-m-2l\}} \notag 
\\ &&\times \binom{k-1}{j-1-m-2l} 
 \frac{(-1)^{j-1+m+l}(j-1)!}{m! l! 2^l}  u^{2m+2l} .  \notag
\eeq
For details of this and other cases see Chapter 15 of \cite{RFG} or Chapter 4 of \cite{ARF}.

\vskip0.6cm \hrule height  0.5pt
\section{The main result - an expectation formula}
\label{sec:expectation}

We now come to the main result of the new theory, the fundamental starting point of 
which was in Jonathan Taylor's 2001 McGill PhD thesis (cf.\ \cite{Taylor-thesis,Taylor1}).

We present it first in the isotropic case. That is, we assume that  $f\:M\to\Rk$
is a  Gaussian random field, the components
of which are independent and identically distributed,  have mean zero and unit variance, 
and are isotropic, with second spectral moment $\lambda_2$ defined by
\beq
\label{lambda2}
\lambda_2 = \E\Big\{\Big(\frac{\partial f(x)}{\partial x_i}\Big)^2\Big\}= 
 \frac{\partial^2 C(x)}{\partial x_i^2}\Big|_{x=0}.
\eeq
(By isotropy, this definition does not depend on $i$.)

Then, for parameter spaces $M$ of dimension $N$, and hitting sets $D$ of dimension $k$,
and under the regularity assumptions of Section \ref{regularity:sec}, we have
\beq
  \label{eq:isotropic:kff:1d}
&&\E \left\{\lips_i\left(A_D(f,M)\right)\right\} \\ &&\qquad\qquad \qquad =
 \sum_{j=0}^{N-i} \sqbinom{i+j}{j} \frac{\lambda_2^{i+j}}{(2\pi)^{j/2}}\lips_{i+j}(M) 
\, \mink_j^{k}(D),
\notag
\eeq
where the so-called combinatorial `flag coefficients'
are defined by
\beqq
\sqbinom{n}{j} = \binom{n}{j} \frac{\omega_n}{\omega_{n-j} \; \omega_j}.
\eeqq

Remaining for the moment in the isotropic case, we shall devote some space to explaining
the meaning of \eqref{eq:isotropic:kff:1d} by considering some examples. 

First of all, we take $k=1$, so that we are dealing with a real valued random field, and
also allow its variance to be $\sigma^2$, which is not necessarily 1. To make life even 
easier, we shall take the parameter space to be the $N$-dimensional rectangle $T$ of
\eqref{rectangles}, and shall concentrate on the case $i=0$ in  \eqref{eq:isotropic:kff:1d},
so that we are looking only at Euler characteristics, rather than general LKCs. 
Finally, we shall take $D=[u,\infty)$. Thus, in summary, we are looking at the expected
Euler characteristic of the excursion set, in a rectangle, of $f$ above the level $u$. 

To substitute into the right hand side of \eqref{eq:isotropic:kff:1d} we need to know how
to  compute the $\lips_{j}(T)$ and $\mink_j^{1}(D)$. The latter, however, are given 
by \eqref{tubes:minkequalshermite:equation} and it is not hard to compute the former from
Steiner's formula. However, to describe them, we need some notation.  

In particular, we define a  {\it face}
 $J$ of $T$, of dimension $k$,  by fixing
$N-k$ coordinates, each of which is fixed at either the top or bottom of the rectangle.
Thus, taking the three dimensional cube as a simple example, it has one three dimensional
`face', this being the full cube itself. It has six two dimensional faces, being its sides.
The twelve edges are one dimensional faces, and the eight corners are zero dimensional faces.
Now let  $\cO_k$ denote the ${N\choose k}$ $k$-dimensional faces which include the origin.
Then it is not to hard to show, from Steiner's formula, that 
\beqq
\lips_j\Big( \prod_{i=1}^N [0,T_i] \Big)=
\sum_{J\in \cO_j} |J|,
\eeqq
where $|J|$ is the ($j$-dimensional) volume of $J$.

Putting the pieces together now gives that, for isotropic random fields,
 $\E\{\p (A_u(f,T)) \}$, or $\E\{\lips_0 (A_u(f,T)) \}$, is given by
\beq
\label{rgeometry:EEC:equation} 
e^{\frac{-u^2}{2\sigma^2}} \sum_{k=1}^{N}
\sum_{J\in \cO_k}
\frac{|J|\,\lambda_2^{k/2}}{(2\pi)^{(k+1)/2}\sigma^{k}}\, 
H_{k-1}\left(\frac{u}{\sigma}\right) 
+\Psi\left(\frac{u}{\sigma}\right), \qquad 
\eeq
or, if $T$ is a cube of side length $T$ rather than a rectangle, by
\beq 
\label{rgeometry:EEC-isotropic:equation} 
e^{\frac{-u^2}{2\sigma^2}}\sum_{k=1}^{N}
\frac{{N\choose k}T^{k}\lambda_2^{k/2} }{(2\pi)^{(k+1)/2}
\sigma^{k}}\, 
H_{k-1}\left(\frac{u}{\sigma}\right)
+\Psi\left(\frac{u}{\sigma}\right).\qquad
\eeq

\begin{figure}[!ht]
{\resizebox{1.5in}{1.5in}{\includegraphics{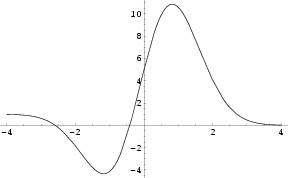}}
\vspace*{-1.25truein}
\hspace*{0.05truein}
\resizebox{1.5in}{1.5in}{\includegraphics{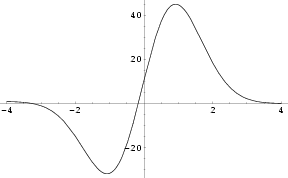}}}
\vspace*{1.65truein}
\caption{Two EEC curves $\E\left\{\p \left(A_u\right)\right\}$ for  Gaussian fields in 
two dimensions\newline for different values of the second spectral moment $\lambda_2$. 
The horizontal axis\newline gives values of $u$ and the vertical axis the EEC.
}
\label{rgeometry:EEC2:figure}
\end{figure}

To get a better feel for this equation, let us look at the cases $N=2$ and $N=3$, taking
$\sigma^2=1$ for simplicity. In the two  dimensionsal case, we obtain
\beq
\label{rgeometry:EEC-isotropic2:equation} 
 \E\left\{\p \left(A_u\right)\right\} =
\left[
\frac{T^2\lambda_2}{(2\pi)^{3/2}}\, u 
\ +\ \frac{2T\lambda_2^{1/2}}{2\pi} 
\right] e^{\frac{-u^2}{2}}  +  \Psi(u).\quad
\eeq
Figure \ref{rgeometry:EEC2:figure} gives two examples, both over the unit square, one with 
$\lambda_2=200$ and one with $\lambda_2=1,000$.

\begin{center}
\begin{tabular}{|c|}   \hline   \\
{\ \ Place Figure \ref{rgeometry:EEC2:figure} near here.\ \  }
\\ \\   \hline \end{tabular}
\end{center}

There are at least four general points that come out of from this example:

\begin{itemize}
\item[(i)] You should note, for later reference, how the expression before
the exponential term can  be thought of as one of a number
of different power series;
one  in $T$, one in $u$, and one in $\sqrt{\lambda_2}$.
\item[(ii)]\ The geometric meaning of the 
negative values of  \eqref{rgeometry:EEC-isotropic2:equation} 
are worth understanding. They are due to the excursion sets having, in the
mean, more holes than connected components for (most) negative values of $u$.
\item[(iii)]\ \ Note the impact of the value of the second spectral moment.
 Large spectral moments lead to 
local fluctuations generating large numbers of small islands (or lakes,
depending on the level at which the excursion  set is taken) and this
leads to larger variation in the values of $\E\{\p (A_u)\}$.
\item[(iv)] \ \ Note that, as $u\to\infty$,  $\E\{\p (A_u)\}\to 0$. This is reasonable,
since at high levels we expect that the excursion set will  be empty, and so have zero Euler
characteristic.

On the other hand as 
$u\to-\infty$, $\E\left\{\p \left(A_u\right)\right\}\to  1$.
The  excursion set geometry behind this is simple.
Once $u<\inf_Tf(t)$ we have
 $A_u\equiv T$, and so $\p \left(A_u\right)=\p (T) =1$.
\end{itemize}

In three dimensions,
\eqref{rgeometry:EEC-isotropic:equation} becomes, again for $\sigma^2=1$,
\beqq
&&
 \E\left\{\p \left(A_u\right)\right\} = \Psi(u)\\ \notag
&&\ \ + e^{\frac{-u^2}{2}}
\left[\frac{T^3\lambda_2^{3/2}}{(2\pi)^{2}}\, u^2 + 
\frac{3T^2\lambda_2}{(2\pi)^{3/2}}\, u 
 + \frac{3T\lambda_2^{1/2}}{2\pi}  - \frac{T^3\lambda_2^{3/2}}{(2\pi)^2} 
\right].
\eeqq
Figure \ref{EEC3-Gauss} gives an example, over the unit cube, with 
$\lambda_2=880$. 
\begin{center}
\begin{tabular}{|c|}   \hline   \\
{\ \ Place Figure \ref{EEC3-Gauss} near here.\ \  }
\\ \\   \hline \end{tabular}
\end{center}

Note that once again there are a number of different power series appearing 
here, although now, as opposed to the two  dimensional case,  there is
no longer a simple correspondence between the powers of $T$, $\sqrt{\lambda_2}$
and $u$.

The  two positive peaks of the curve are due to $A_u$ being primarily
composed of a  number of 
simply connected components for large $u$ and primarily of simple 
holes for negative $u$.
(Recall that in the three dimensional case the Euler characteristic of a set
is given by the number of components minus the number of handles plus the 
number of holes.) The negative values of $\E\{\p (A_u)\}$ for $u$ near zero
are due to the fact that $A_u$, at those levels, is composed mainly of a 
number of interconnected, tubular like regions; i.e.\ of handles.

Although we have discussed only the real valued case in detail above, and for rather simple 
hitting sets, the extension of the isotropic case to vector valued and to Gaussian related
random fields should now be clear. Essentially, the argument is the same, but the 
computation of the GMFs $\cal M_j^k$ becomes more complicated. Consider the case of a
real valued $\chi^2$ random field, given by
\beqq
f(x)=g_1^2(x)+\dots +g_k^2(x),
\eeqq
where the $g_j$ are independent, identically distributed, 
isotropic, mean zero, Gaussian  random fields with second spectral moment $\lambda_2$.
Then the expectations $\E\{\lips_i(A_u(f,M))\}$ follow immediately from 
\eqref{eq:isotropic:kff:1d}. The $\mink_j^{k}$ are as given in  \eqref{chisqared}, and the 
$\lips_j$ are precisely as we just calculated for the simple Gaussian case.

An example of what the corresponding expected Euler characteristic curve looks like is 
given in Figure \ref{EEC3-chi}, for which  we have taken 
$k=5$ and $\lambda_2$ (the second spectral moment of the underlying Gaussian processes)
to be 20. (The reason for this particular choice of parameters will be given in 
Section \ref{sec:examples}.) 

\begin{center}
\begin{tabular}{|c|}   \hline   \\
{\ \ Place Figure \ref{EEC3-chi} near here.\ \  }
\\ \\   \hline \end{tabular}
\end{center}

The $\chi^2$ example exemplifies one of the most important properties of 
\eqref{eq:isotropic:kff:1d}, which is a  {\em separation of parameters}, by which we mean
that the parameters relating to the random field, to the parameter set, and to the 
hitting set do not interact in other than a multiplicative fashion.

Unfortunately, this is not quite true when we leave the isotropic situation. In the general
case -- with assumptions as above but with isotropy (and the implied stationarity) replaced
by an assumption only of constant unit variance -- \eqref{eq:isotropic:kff:1d} becomes
\beq
  \label{eq:gaussian:kff:1d}
 &&\E \left\{\lips_i\left(A_D(f,M)\right)\right\}\\ && \qquad\qquad \qquad  =
 \sum_{j=0}^{N-i} \sqbinom{i+j}{j} (2\pi)^{-j/2}\lips_{i+j}(M) 
\, \mink_j^{k}(D).
\notag
\eeq
However, in this equation the LKCs no longer have the simple meanings that they have had
up until now. We shall try, in one short paragraph, to explain how they change, but the details
are technical and so we refer the interested reader to \cite{RFG} for details.

Consider either  a real valued Gaussian random field or a single component of a vector
valued field. In either case, denote this by $f$. At a point $x\in M$, take two unit vectors,
say $V_1$ and $V_2$, with base at $x$. Let $V_if$ denote the directional derivative of 
$f$ in direction $V_i$. Now define a function $g_x$ on all such vectors by setting
\beq
\label{metric:equn}
g_x\left(V_1,V_2\right) =\E\left\{ V_1f(x)\,V_2f(x)\right\}.
\eeq
Then $g_x$ defines what is known as a Riemannian metric on $M$, turning it into a Riemannian
manifold, and with a corresponding notion of volume. It turns out that tube formulae such as
\eqref{steiner} still hold for such volumes, and the coefficients in the corresponding
expansion are the new LKCs. Note that the LKCs (except for the Euler characteristic
$\lips_0$) now incorporate information about the 
random field $f$, specifically information about the covariances of partial derivatives.
The fact that the Euler characteristic does not change is a deep result known as the 
Gauss-Bonnet theorem, and, as we shall see in the following section, is extremely
convenient when looking at excursion probabilities.

Despite the entrance of Riemannian geometry into our story, there are at least two cases in
which life is not too bad, and it is possible to compute the LKCs quite simply.
For all the other cases we shall describe, in Section \ref{estimating:section},
 how to estimate them from data. In fact, this estimation procedure is so simple that
it enables one to use formulae like  \eqref{eq:gaussian:kff:1d}, at least for the 
important case $i=0$, without even having to know the definition of the  LKCs appearing on
the right hand side of the equation. 

Now, however, we consider the two simple cases.
 The first occurs under stationarity, for which we can define
a collection of second order spectral moments by
\beq
\label{lambdaij}
\lambda_{ij} = \E\left\{\frac{\partial f(x)}{\partial x_i} \frac{\partial f(x)}{\partial x_j}
\right\}= 
- \frac{\partial^2 C(x)}{\partial x_i\partial x_j}\Big|_{x=0}.
\eeq
Return now to the case \eqref{rectangles} of a $N$-dimensional rectangle $T$.
Let $(e_1,\dots,e_N)$ be the usual set of axes in $\RN$ and for a $k$-dimensional face  
$J$ let $\Lambda_J$ be the matrix of second spectral moments $\lambda_{ij}$ for which
$i$ and $j$ satisfy the requirement that both $e_i$ and $e_j$ intersect $J$ at  
points other than the origin.
(If $\dim (J)=0$, then we define $\Lambda_J$ to be the constant 1.)  Then the LKCs of
\eqref{eq:gaussian:kff:1d} become
\beqq 
\sum_{J\in \cO_k}
|J|\,\left[\det\left(\Lambda_J\right)\right]^{1/2},
\eeqq
something which is quite simple to evaluate in practice. 
Unfortunately, however, the computation is a little more complicated when the parameter
space is not a rectangle. However, the GMFs, as always, are independent of the properties
of both $f$ and $M$.

Another helpful situation arises when one is interested in high level excursion sets of 
real valued random fields. Looking back at results such as 
\eqref{rgeometry:EEC:equation} and \eqref{rgeometry:EEC-isotropic:equation}, we note
that they have the form of a power series expansion in $u$, so that, if $u$ is large, 
the first term in the expansion should dominate all the others. Putting this together with 
\eqref{eq:gaussian:kff:1d} now tells us that, at least for real valued Gaussian fields of
unit variance, we should have 
\beq
  \label{LK-asymp}
&& \E \left\{\lips_i\left(A_u(f,M)\right)\right\}
\\ && \qquad 
\approx \sqbinom{N-i+j}{N-i} (2\pi)^{-(N-i)/2}\lips_{N}(M) \, \mink_{N-i}^{1}([u,\infty)),
\notag
\eeq
where we intepret $a(u)\approx b(u)$ to mean that $\lim_{u\to\infty} a(u)/b(u)=1$.
 Thus, all we need, at least for asymptotics, is to know how to compute
$\lips_N(M)$. However, although the other LKCs are hard to compute in general, there is a
simple representation for this one. In fact, if we let $\Lambda(x)$ be the $N\times N$
matrix with elements
\beqq
\lambda_{ij}(x) = \E\left\{\frac{\partial f(x)}{\partial x_i} \frac{\partial f(x)}{\partial x_j}
\right\}= 
 \frac{\partial^2 C(u,v)}{\partial u_i\partial v_j}\Big|_{u=v=x},
\eeqq
then one can show that 
\beqq
\lips_N(M) = \int_M \left[\det\left( \Lambda(x) \right)\right]^{1/2}\, dx,
\eeqq
so that \eqref{LK-asymp} becomes a very general, and useful, result.

\vskip0.6cm \hrule height  0.5pt
\section{Excursion probabilities} 

We now want to see how the formulae of the previous section 
can be applied to evaluating excursion probabilities, something which we have already
mentioned is one of their main applications.

Actually, most of what we had to say was already said in equation \eqref{ECequivalence}, where
we pointed out that, for many random fields,  
\beq
\label{ECequivalence1}
&& \Big|\P\Big\{\sup_{x\in M} f(x) \geq u\Big\} - 
\E\left\{\p \left(A_u(f,M)\right)\right\}\Big|
\leq error(u),
\notag \\
\eeq
where $\p =\lips_0$ is the Euler characteristic.
The previous section was devoted to ways of computing the expectation here for a wide class
of random fields, and so all we have left to do is to say something about the error function and
offer a heuristic explanation as to why such a result might hold.
We start with the heuristics.

Note firstly that 
\beq
\label{equivalence}
\sup_{x\in M} f(x)\geq  u &\iff&  A_u(f,M)\neq \emptyset\\
&\iff& \#\{\text{connected components of $A_u$}\} \geq 1. \notag
\eeq
Consider the Euler characteristic of $A_u$. Recall, that we already saw, in two and 
three dimensions, that $\p (A_u)$ is equal to the number of connected components of $A_u$, plus or
minus other factors,  such as the numbers of holes and handles. In fact, it is true in any
dimension that the Euler characteristic of a set is given by the number of its connected 
components, plus an alternating sum of what might be called `higher order' functionals, which,
as the dimension grows, describe more and more complicated geometric structures. However,
it can be shown that, as $u$ becomes large, the structure of an excursion set $A_u$ becomes rather
simple, containing only simply connected components, without holes, handles, etc. In fact, if
$u$ is large enough, $A_u$ will contain only one (small) simply connected component, and so 
will have Euler characteristic of one. In other words, for large $u$ we can replace 
\eqref{equivalence} by
\beq
\label{as-equivalence}
\sup_{x\in M} f(x)\geq  u \ \   \mbox{``$\iff$''} \ \  \p\left(A_u(f,M)\right) =1,
\eeq
where we interpret ``$\iff$''
to mean `asymptotically if and only if', and
now $\p(A_u)$ can, with high probability, only  take the values 0 or 1.

Taking probabilities on both sides of this equation, and recalling that if $X$ is
a 0-1 random variable then $\P\{X=1\}=\E\{X\}$, gives that 
\beqq
\P\Big\{\sup_{x\in M} f(x) \geq u\Big\} \approx 
\E\left\{\p \left(A_u(f,M)\right)\right\}
\eeqq
for large $u$, or, equivalently, that the error term in 
\eqref{ECequivalence1} should be small for large $u$.

This heuristic argument works well for most random fields, and you can read more about it and
other heuristic arguments for random field extrema in \cite{ADLREV}. However,
 a rigorous result has been proven  only for
constant variance Gaussian fields. (Unfortunately, though, the above heuristics have little to
do with the proof.) We shall give just one example of it.

Suppose $f$ is Gaussian, isotropic, and with constant unit variance. Then the error function
in  \eqref{ECequivalence1} is of the form
\beq
\label{errorf}
error(u)=e^{-{u^2} (1+ \sigma^{-2}_c )/2},
\eeq
where
\beqq
\sigma^2_c =
\Var\left(\frac{\partial^2 f(x)}{\partial x_1^2}\Big| f(x)\right) =
 \frac{\partial^4C(x)}{\partial x_1^4}\Big|_{x=0}\, -\, 1.
\eeqq
This particular result is an example of a very general phenomenon for 
unit variance Gaussian random fields, of basically the
same form, but in which the exponent on the right hand side of  \eqref{errorf} is
replaced by $- u^2(1+\eta)/2$, for a $\eta>0$ that is not always easy to identify. 
However, even when $\eta$ cannot be identified, this is a remarkable result.

To see why, consider, for example, the unit variance stationary case with parameter space a $N$-dimensional
rectangle (although any set would do)  in which case we now know from 
\eqref{rgeometry:EEC:equation} that we can rewrite the above as
\beqq
\P \Big\{\sup_{x\in T} f(x) \geq  u \Big\}
=  C_0\Psi\left(u\right)  +
u^{N} \, e^{-u^2/2}
\sum_{j=1}^{N}     C_j u^{-j},   
\eeqq
plus an error of order 
$ o(e^{- u^2(1+\eta)/2})$
where the $C_j$ are constants depending on the LKCs of $T$ and the second order
spectral moments of $f$.

Why is this remarkable? If we think of the right hand side above  as
a terminated version of an infinite power series expansion, 
then it would be natural to expect that the error term  would be of the order 
of the `next' term in the expansion,  and so of order 
$ u^{-1} e^{-u^2/2}$. However, we have seen that 
this is {\it not} the case, and that the error is actually 
{\it exponentially smaller} than this. Presumably this theoretical result
is behind the fact, well established in practice, that the approximation of
excursion probabilities by expected Euler characteristics of excursion sets
works remarkably well. In general,  numerical and other evidence suggests 
that at levels at which the
true excursion probability is less than 0.10, the error in the Euler characteristic
approximation is no more than 5-10\% of the true probability.

\vskip0.6cm \hrule height  0.5pt
\section{On measuring, calculating and estimating LKCs and Minkowski functionals}
\label{estimating:section}
 
We now know that there are three main components that appear in the new theory that 
has been discussed above. The first are the \LKC s of excursion sets, which are random 
variables and can be computed from data. We have not gone into detail as to how to do this,
and it is not always easy. However, it {\it is} easy for the 
Euler characteristic $\p$ ($\equiv \lips_0$)  and we shall see now how to do this.

Suppose, which is almost always that case,  that we see our data not on a continuum,
but rather on a rectangular lattice $L_\delta$, of edge size $\delta$, so that 
rather than seeing $\{f(x)\}_{x\in M}$ we see only $\{f(x)\}_{x\in L_\delta\cap M}$.
Then we also cannot really see the true excursion sets $A_D(f,M)$, but rather 
rectanglar approximations to them of the form
\beqq
A_D^\delta(f,M) \definedas \left\{x\in M\cap L_\delta\:f(x)\in D \right\}.
\eeqq
If we join each point in  $A_D^\delta(f,M)$ to those distance $\delta$ from it (there
are a maximum of $2N$ such points) we
construct an $N$-dimensional object, made up of $N$-dimensional squares, 
$(N-1)$-dimensional faces, and so on, down to 0-dimensional vertices (which are actually
the points of $A_D^\delta(f,M)$. Denoting the number of such $k$-dimensional faces by
$N_k^\delta$, $k=N,N-1,\dots,0$, it turns out that 
\beqq
\p\big( A_D^\delta(f,M)\big) =\sum_{k=0}^N (-1)^{N-k}   N^\delta_k.
\eeqq
Not surprisingly, under mild regularity conditions, it also turns out that 
$\p( A_D^\delta(f,M)) \to \p( A_D(f,M))$ as $\delta\to 0$. This is how one computes,
or approximates, the Euler characteristic of excursion sets in practice. The 
remaining LKCs, particularly in the non-isotropic case, are somewhat more complex.
In the isotropic case, however, one can use the definitions given in the Introduction
to compute them.

The remaining two components of the theory are the LKCs of the parameter space, 
measured with respect to the Riemannian metric \eqref{metric:equn} as described at the
end of Section \ref{sec:expectation}, 
and the Gaussian Minkowski functionals of the hitting set $D$, described in 
Section \ref{gmf:sec}. Computing either of these can involve considerable knowledge of 
differential geometry. The good news is that there are ways around this.

First of all, the  Gaussian Minkowski functionals depend only on the hitting set $D$,
and there are only a limited number of hitting sets that are in common usage. Thus
one can generally find what one needs already done somewhere. An exhaustive list of 
known examples can be found in Chapter 4 of \cite{ARF}, along with techniques for 
developing any additional cases that might arise.

The case of the \LKC s is somewhat different. If the random field is stationary, and 
the parameter space $T$ is a rectangle, then we are covered by 
\eqref{rgeometry:EEC:equation}. A few other examples are given in Chapter 4 of 
\cite{ARF}. All of these examples require estimation of the parameters $\lambda_{ij}$,
or, at least, of determinants of matrices made up by them. A more direct way of 
estimating LKCs has been proposed in \cite{Taylor-Worsley-JASA}. The strength of this
approach is that there is actually no need for the practitioner to have any 
serious depth of understanding of the geometrical meaning of the $\lips_j$. Rather,
they can be treated almost as nuisance parameters to be estimated from data via a rather
straightforward algorithm.

\vskip0.6cm \hrule height  0.5pt
\section{An application: Model identification}
\label{sec:examples} 

To conclude this review, I want to describe one rather interesting and powerful
 application of the new formulae of Section \ref{sec:expectation}, and rather
than describing the general problem, I shall do it via a particular example.

\subparagraph{The example}

Suppose that we are given a real valued random field $f$ on $\real^3$, which we observe on
the unit cube $[0,1]^3$. Our problem is that we are not certain what the distribution of
the field is. For simplicity, we shall assume that the choices are Gaussian and $\chi^2$,
where both are assumed to be stationary and isotropic.

Then one way to choose between the two is to calculate, from the data, an empirical
Euler characteristic curve, given by
\beqq
\widehat \p(u)  = \p \left(A_u\left(f,[0,1]^3\right)\right). 
\eeqq

\begin{center}
\begin{tabular}{|c|}   \hline   \\
{\ \ Place Figure \ref{EEC3-chi-Gauss} near here.\ \  }
\\ \\   \hline \end{tabular}
\end{center}

\begin{figure}[!ht]
{\resizebox{3.7in}{1.7in}{\includegraphics{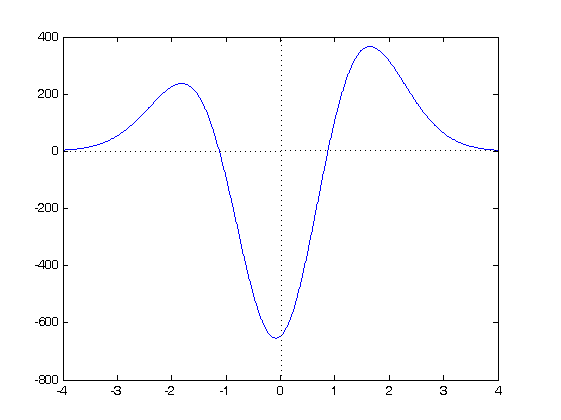}}}
\caption{EEC for a Gaussian field in three dimensions. 
 Axes as for Figure \ref{rgeometry:EEC2:figure}.}
\label{EEC3-Gauss}
\end{figure}

If the data is Gaussian, then the graph of $\widehat \p$ should be of the general
form of Figure \ref{EEC3-Gauss}, where the actual height and width of the
graph will depend on the
variance and second spectral moment of $f$.

\begin{figure}[!ht]
{\resizebox{3.7in}{1.7in}{\includegraphics{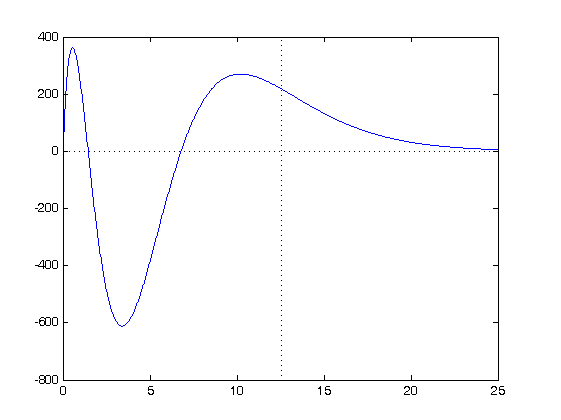}}}
\caption{EEC for a $\chi^2_5$ field in three dimensions.
Axes as for Figure \ref{rgeometry:EEC2:figure}.}
\label{EEC3-chi}
\end{figure}

If the data is $\chi^2$, then the graph of $\widehat \p$ will have a different form,
closer to that of Figure \ref{EEC3-chi}, for example, which recall is the expected
Euler characteristic graph for a $\chi^2$ random field with 5 degrees of freedom, and
with variance and spectral moment chosen to match those of the Gaussian random 
field of   Figure \ref{EEC3-Gauss}.

The graphs of  Figure \ref{EEC3-Gauss} and  Figure \ref{EEC3-chi} are quite different,
and thus one should be able to see the same qualitative differences in graphs of 
$\widehat \p$ in the two cases. This indeed is generally the case, and so we see
that the expected Euler characteristic (EEC) curve provides a diagnostic tool for differentiating between models.

One caveat should be made here, however, before we continue. The significant differences
between the Gaussian and $\chi^2$ graphs would, of course, weaken considerably were
we to take the latter with a far higher degree of freedom, say 30 or more. This,
however, is far from surprising, since in such a regime the central limit theorem tells
us that the not only are the differences between Gaussian and $\chi^2$ random variables
rather small, the same is true for  Gaussian and $\chi^2$ random fields.

Returning to the $\chi^2_5$ example, suppose we were to  adopt what we claimed at the
beginning of Section \ref{sec:rfields} was the erroneous approach of `transforming to
normality' via a transformation as in  \eqref{intro:standard:equn}. In fact,
let us go a little further now, assuming that we know our data is sampled from 
a $\chi^2_5$ random field $f$. Let  $F_5$ to be the (known)  
distribution function of a $\chi^2_5$ random variable, and create a new random field
by setting
\beq
\widetilde f(x) \definedby \Phi^{-1}\big(\widehat F_5(f(x))\big).
\eeq
In other words, we perform the `Gaussianisation' of \eqref{intro:standard:equn} which,
because of our assumptions, we can do without data.

Then the univariate distributions of $\widetilde f$ are standard Gaussian, and so it
would be impossible, looking merely at these, to distinguish between $\widetilde f$ and a 
Gaussian random field.

\begin{figure}[!ht]
{\resizebox{3.7in}{1.7in}{\includegraphics{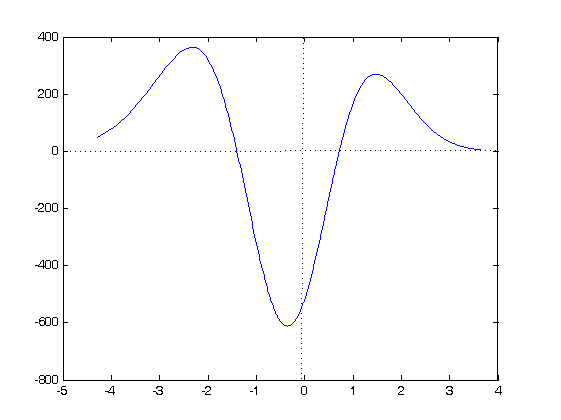}}}
\caption{EEC for a `normalised' $\chi^2_5$ field in three dimensions.
Axes as for Figure  \ref{rgeometry:EEC2:figure}.}
\label{EEC3-chi-Gauss}
\end{figure}

Figure  \ref{EEC3-chi-Gauss}  gives  the EEC curve for  $\widetilde f$, with parameters
chosen so that $\widetilde f$ has the same variance and second spectral moment of the
Gaussian field of Figure \ref{EEC3-Gauss}. This curve is quite different to the
corresponding Gaussian one, both in terms of its general shape and the rate  at which it
reaches its asymptotes. Thus, while univariate marginal distributions of the 
Gaussianised random field offer no indication of inherent non-Gaussianity, the 
EEC curve does.

This fact has often been used, particularly in the astrophysical literature, not only
as a way to distinguish between two models, but to search among models for the 
correct one (e.g.\ \cite{Sloan,TORRES,VPGHG}).

\subparagraph{What we can learn from the example}

After reading the above `application', a referee asked a number of pertinent questions
that would probably be of interest to any reader who works with spatially correlated
data. Thus, it seems appropriate to close with a paraphrased, concatenated
 version of his questions,
along with the best answers that I can provide to them.

{\it 1: What are the advantages of the above inferential procedure to one based on standard 
distributional tests?}

In a certain sense, a test for normality based on Euler characteristics or, indeed,
any of the LKCs, {\it is} a distributional test. The basic expectation
formulae depend only the variance of $f$ and the variances and covariances of its
first order partial derivatives. If, for example, the random field is assumed
to be isotropic, as in the example above, then the test described there is really
only a test of the joint normality or otherwise of the distribution of the
Gaussianised field and its first order derivatives.

What makes it somewhat different from standard tests is {\it how} it tests for 
non-normality, by placing emphasis on the geometric structures that are generated by
the random field. Whether or not this is a good thing will depend very much on where
one's interest in the data is concentrated. If it is concentrated on geometry, we 
believe that EEC techniques provide a more powerful approach to testing for non-normality
that standard tests. This belief certainly seems to be shared, and put into practice,
 by large parts of the  brain mapping and astrophysics communities for whom
the analysis of  shape has become an important part of their inferential procedures.

Note, however, that if assumptions of neither isotropy nor stationarity are made, then 
straightforward distributional tests are impossible to apply, since there are
an uncountable infinity of different distributions to be considered. In these cases
a summary statistic of some sort must be employed. The geometric approach suggested 
above seems to be particularly suited to these scenarios.

{\it 2: How does the theory presented here convey information about the long range  properties of the random field?}

I have left this point to last, since it is an extremely important one. The fact is
that none of the formulae of this paper are affected by, or descriptive of,
long range properties of the random field. Of course, all the formulae for the 
means $\E\{\lips_j(A_u(M,f)\}$ depend on the size and shape of $M$, and, for example,
 they grow,
in an additive fashion, as $M$ grows. However, this growth does not take into account
the rate of decay, for example, of the covariance function of $f$.

At first sight, this is rather surprising, but, in fact, it is not something particularly
new. Consider, for example, the excursion probability \eqref{excursionprob}. From the 
asymptotic equivalence \eqref{ECequivalence} (or any other technique for computing 
\eqref{excursionprob})
and what we now know about the mean Euler characteristic of excursions sets,
it follows that, at least a high levels, the excursion probability also 
does not depend on the long range behaviour of the 
covariance function. There are many other examples in the theory of Gaussian random
fields for which this kind of behaviour occurs, and the geometric theory of this paper
is just one more.

Is this a strength or a weakness? The answer, of course, depends on what one is interested
in. If long range behaviour of covariance functions is important for some other
reason, then the theory of this paper has little to offer to help in understanding
or estimating them. On the other hand, it is nice to know that there are many  
interesting aspects of random fields that do not depend on the decay rates of 
covariance functions. Indeed, it 
may well be rather comforting, since such decay rates are one of the hardest 
phenomena to accurately measure, in practice, from real data.

\vskip0.6cm \hrule height  0.5pt

\section{Acknowledgements and apologies} Throughout this review I have given few
references and very little detail as to the long history of this subject. Two names
stand out. One is Jonathan Taylor, whose fundamental contributions to the theory have
 already been mentioned, and 
the second is Keith Worsley, who has made more contributions
to the application of the Euler characteristic approach to random fields than anyone
else I know. However, there are {\it many} other names that deserved to be 
mentioned, and were not.
I acknowledge that I owe an apology to both these authors and the reader for not mentioning
them by name.

My justification for this oversight is two-fold. Firstly, 
in a brief review article like this one it would be
impossible to even attempt to do credit to the history of the subject. Secondly, there
are sources that you can turn to for the missing details, including 
\cite{ADLREV,RFG} and \cite{ARF}.

Finally, a debt of gratitude is due to Einat Engelhart, a Technion MSc student, who
helped  with some of the numerical work.

\vskip0.6cm \hrule height  0.5pt

\end{multicols}

\newpage

\begin{thebibliography}{99.}

\bibitem{ADLREV}   Adler, R.J.  On excursion sets, tube formulae, and maxima of random
fields,  {\it
Annals of Applied Prob}, {\bf 10}, (2000), 1--74.

\bibitem{AMR}   Adler, R.J.,  Muller, P.,  Rozovskii, B.  (eds.) (1996) {\it Stochastic
Modelling in
Physical Oceanography} (1996), Birkha\"user, Boston.


\bibitem{RFG}
 Adler, R.J.  and Taylor,  J.E.   (2007) {\em Random Fields and Geometry.} Springer, Boston.

\bibitem{ARF}
 Adler, R.J., Taylor, J.E.   and  Worsley, K.J.  (2009?) {\em Applications of 
Random Fields and Geometry: Foundations and Case Studies}, 
In preparation.  Some chapters already avaliable at  ie.technion.ac.il/Adler.phtml.


\bibitem{Christakos1}  Christakos, G., (1991) {\em Random Field Models in Earth Sciences},
 Academic.

\bibitem{Christakos2} 
 Christakos, G., (2000) {\em Modern Spatiotemporal Geostatistics},. Oxford Univ.\ Press,
Oxford.

\bibitem{Daley} Daley, R., (1991) {\em Atmospheric Data Analysis}, Cambridge Univ.\ Press,
Cambridge.

\bibitem{pet} 
{Evans, A. C. and  Marret, S. and  Neelin, P. and
Collins, L. and Worsley, K.J. and  Dai, W. and Milot, S. and Meyer, E. and
and Bub, D.} (1992) {Anatomical mapping of functional
activation in stereotactic coordinate space}, 
{\em Neuro{I}mage},
{1}, 
{43--53}.

\bibitem{Sloan}
Gott, J.R., Hambrick, C.D., Vogeley, M.S., Kim, J.,
Park, C., Choi, Y-Y., Cen, R., Ostriker, J.P., and Nagamine, K.
(2008)
Genus topology of structure in the Sloan digital sky survey: Model testing
{\it The Astrophysical Journal}, 675, 16--28.



\bibitem{MATERN}
{Mat{\'e}rn, B.}, (1960)   {\em Spatial variation: {S}tochastic models and their application
              to some problems in forest surveys and other sampling
              investigations}, {Meddelanden Fran Statens Skogsforskningsinstitut, Band 49, Nr.
              5, Stockholm}.

\bibitem{Olea} Olea, R. A., (1999) {\em Geostatistics for Engineers and Earth Scientists},
Kluwer Acad.\ Publ., Boston, MA.

\bibitem{Rubin} Rubin, Y., (2002) {\em Applied Stochastic Hydrogeology}, Oxford Univ.\ Press, 
New York.

\bibitem{SampsonGuttorp}
Sampson, P.D. and Guttorp, P.
 {Nonparametric estimation of nonstationary spatial covariance structure},
 {\it J.\ Amer.\ Stat.\ Assoc.}
      {\bf 87},
       {1992},
      {108--119}.
  

\bibitem{SSSW}  Shafie,  K., Sigal, B., Siegmund, D.  and  Worsley, K.J. (2003)
	 {Rotation space random fields with an application to f{MRI} data},
 {\em Ann. Statist.}, {31},	 {1732--1771}.
  
\bibitem{Smootbook}  Smoot, G. and Davidson, K., (1993)
     {\em Wrinkles in {T}ime},
{William Morrow},
{New York}

\bibitem{Taylor-thesis} 
Taylor, J.E. {\it Euler Characteristics for
    Gaussian Fields on Manifolds}, Ph.D.\ dissertation,  McGill, (2001).
%

\bibitem{TTA} {Taylor, J.E., Takemura, A. and Adler, R.J.} (2005)
 {Validity of the expected {E}uler characteristic heuristic},
{\em Ann.\ Probab.},  {33}, {1362--1396}.
 

\bibitem{Taylor-Worsley-JASA}
 {Taylor, J.E. and Worsley, K.J.},
 {Detecting sparse signals in random fields, with an application
              to brain mapping},
 {\em J.\ Amer\. Statist.\ Assoc.},
 {2007},
 {479},
 {913--928}.



\bibitem{Taylor1}
Taylor, J.E. (2006)  A Gaussian kinematic formula, {\em Ann. Probability}, 34, 122--158.


\bibitem{TORRES}
{Torres, S.} (1994)
{Topological analysis of {COBE-DMR} cosmic microwave background maps},
{\em Astrophysical Journal},
    {423}, {L9--L12}.

\bibitem{VPGHG}  
{Vogeley, M.S. and Park, C. and Geller, M.J. and Huchin, J.P. and Gott, J.R.} (1994)
{Topological analysis of the {C}f{A} redshift survey},
{\em Astrophysical Journal},{420}, {525--544}.


\end{thebibliography}
\end{document}